\documentclass[english]{article}
\usepackage[latin1]{inputenc}
\usepackage{amsmath,amsthm,amssymb,babel}

\newtheorem{teo}{Theorem}

\newtheorem{lemma}{Lemma}

\def\proof{{\it Proof.}\ }
\def\endproof{\hfill $\Box$\par\vskip3mm}

\def\eq#1{(\ref{#1})}

\def\neweq#1{\begin{equation}\label{#1}}
\def\endeq{\end{equation}}

\def\phi{\varphi}
\def\RR{{\mathbb R} }

\def\di{\displaystyle}
\def\ri{\rightarrow}

\date{}

\title{\sc On a non-homogeneous eigenvalue problem involving a
potential: an Orlicz-Sobolev space setting \thanks{
Correspondence address: Vicen\c{t}iu R\u{a}dulescu, Department of
Mathematics, University of Craiova,  200585 Craiova, Romania. E-mail:
{\tt
vicentiu.radulescu@math.cnrs.fr}}}
\author{\sc Mihai Mih\u ailescu$\,^{a,b}$ \qquad Vicen\c{t}iu R\u{a}dulescu$\,^{a,c}$ \qquad Du\u san Repov\u s$\,^{d}$\\
\small
$^a\,$Department of Mathematics, University of Craiova,  200585 Craiova,
Romania\\
\small $^b\,$Department of Mathematics, Central European
University,   Budapest,
Hungary 1051\\
\small $^c\,$Institute of Mathematics ``Simion Stoilow" of the Romanian Academy,\\
\small P.O. Box 1-764,  Bucharest, Romania 014700\\
\small $^d\,$Faculty of Mathematics and
Physics, University of Ljubljana, Jadranska 19, Ljubljana, Slovenia 1000\\
\small E-mail addresses: {\tt mmihailes@yahoo.com}\qquad {\tt vicentiu.radulescu@math.cnrs.fr}\qquad {\tt dusan.repovs@uni-lj.si}}

\begin{document}
\baselineskip16pt \maketitle \noindent{\small{\sc Abstract}. In
this paper we study a non-homogeneous eigenvalue problem involving variable
growth conditions and a potential $V$. The problem is analyzed in the context of Orlicz-Sobolev spaces. 
Connected with this problem we also study the optimization problem for the particular eigenvalue given by the
infimum of the Rayleigh quotient associated to the problem with respect to the potential $V$ when $V$ lies in 
a bounded, closed  and convex subset of a certain variable exponent Lebesgue space.
}

\section{Introduction and preliminary results}
Let $\Omega$ be a bounded domain in $\RR^N$ ($N\geq 3$) with smooth boundary
$\partial\Omega$. Assume that  $a_i:(0,\infty)\ri\RR$, $i=1,2$, are two functions such that the mappings $\phi_i:\RR\rightarrow\RR$,
$i=1,2$, defined by
$$\phi_i(t)=\left\{\begin{array}{lll}
a_i(|t|)t, &\mbox{for}&
t\neq 0\\
0, &\mbox{for}& t=0\,,
\end{array}\right.$$
are  odd, increasing homeomorphisms from $\RR$ onto $\RR$, $\lambda$ is a  real number, $V(x)$ is a potential 
and $q_1$, $q_2$, $m:\overline\Omega\rightarrow(1,\infty)$ are continuous functions. We analyze the eigenvalue problem
\begin{equation}\label{1}
\left\{\begin{array}{lll}
-{\rm div}((a_1(|\nabla u|)+a_2(|\nabla u|))\nabla u)+V(x)|u|^{m(x)-2}u=\lambda(|u|^{q_1(x)-2}+|u|^{q_2(x)-2})u, &\mbox{if}&
x\in\Omega\\
u=0, &\mbox{if}& x\in\partial\Omega\,.
\end{array}\right.
\end{equation}
The interest in analyzing this kind of problems is motivated by some recent advances in the study of
eigenvalue problems involving non-homogeneous operators in the divergence form. We refer especially to the results in
\cite{FZZ,mihradproc,manuscripta,F,AA,fourier,blms}. Problem \eq{1} can be placed in the context of the above results since
in the particular case when $q_1(x)=q_2(x)=q(x)$ for any $x\in\overline\Omega$ and $V\equiv 0$ in $\Omega$
it was studied  in \cite{AA}. The form of problem \eq{1} becomes a natural extension of the problem studied in \cite{AA} with the presence of
the potential $V$ in the left-hand side of the equation and by considering that in the right-hand side we can have $q_1\neq q_2$ on
$\overline\Omega$. 

In order to go further we introduce the functional space setting where problem \eq{1} will be discussed. In this context 
we notice that the operator in the divergence form is not homogeneous and thus, we introduce an Orlicz-Sobolev space setting 
for problems of this type. On the other hand, the presence of the continuous functions $m$, $q_1$ and $q_2$ as exponents 
appeals to a suitable variable exponent Lebesgue space setting. In the following, we give a brief description of the Orlicz-Sobolev 
spaces and of the variable exponent Lebesgue spaces.

We start by recalling some basic facts about Orlicz spaces. For more
details we refer to the books by D. R. Adams and L. L. Hedberg
\cite{AHed}, R. Adams \cite{A} and M. M. Rao and Z. D. Ren
\cite{rao} and the papers by Ph. Cl\'ement et al. \cite{Clem1,
Clem2}, M. Garci\'a-Huidobro et al. \cite{Gar} and J. P. Gossez
\cite{G}.

For $\phi_i:\RR\rightarrow\RR$, $i=1,2$, which are odd, increasing homeomorphisms from $\RR$ onto $\RR$, we define
$$\Phi_i(t)=\int_0^t\phi_i(s)\;ds,\;\;\;(\Phi_i)^\star(t)=\int_0^t(\phi_i)^{-1}(s)\;ds,\qquad \mbox{for all}\ t\in\RR,\;i=1,2\,.$$
We observe that $\Phi_i$, $i=1,2$, are {\it Young functions}, that is, $\Phi_i (0)=0$, $\Phi_i$ are convex, and $\lim_{x\ri\infty}\Phi_i (x)=+\infty$.
Furthermore, since $\Phi_i (x)=0$ if and only if $x=0$, $\lim_{x\ri 0}\Phi_i (x)/x=0$, and $\lim_{x\ri \infty}\Phi_i (x)/x=+\infty$, then
$\Phi_i$ are called  {\it $N$-functions}.
The functions $(\Phi_i)^\star$, $i=1,2$, are called the {\it complementary} functions of $\Phi_i$, $i=1,2$, and they satisfy
$$(\Phi_i)^\star (t)=\sup\{st-\Phi_i (s);\ s\geq 0\},\qquad\mbox{for all $t\geq 0$}\,.$$
 We also observe that $(\Phi_i)^\star$, $i=1,2$, are also  $N$-functions and Young's inequality holds true
 $$st\leq\Phi_i (s)+(\Phi_i)^\star (t),\qquad\mbox{for all $s,t\geq 0$}\,.$$

The Orlicz spaces $L_{\Phi_i}(\Omega)$, $i=1,2$, defined by the $N$-functions $\Phi_i$
(see \cite{AHed,A,Clem1}) are the spaces of measurable functions $u:\Omega\ri\RR$ such that
$$\|u\|_{L_{\Phi_i}}:=\sup\left\{\int_\Omega uv\;dx;\ \int_\Omega(\Phi_i)^\star (|g|)\;dx\leq 1\right\}<\infty\,.$$
Then $(L_{\Phi_i}(\Omega),\|\,\cdot\, \|_{L_{\Phi_i}} )$, $i=1,2$,
are Banach spaces whose norm is equivalent to the
 Luxemburg norm
$$\|u\|_{\Phi_i} :=\inf\left\{k>0;\ \int_\Omega\Phi_i\left(\frac{u(x)}{k}
\right)\;dx\leq 1\right\}.$$ For Orlicz spaces  H\"older's
inequality reads as follows (see \cite[Inequality 4, p. 79]{rao}):
$$\int_\Omega uvdx\leq 2\,\|u\|_{L_{\Phi_i}}\, \|v\|_{L_{(\Phi_i)^\star}}\qquad\mbox{for all $u\in L_{\Phi_i}(\Omega)$
and $v\in L_{(\Phi_i)^\star}(\Omega)$},\;i=1,2\,.$$

Next, we introduce the Orlicz-Sobolev spaces. We denote by $W^1L_{\Phi_i}(\Omega)$, $i=1,2$, the Orlicz-Sobolev spaces defined
by
$$W^1L_{\Phi_i}(\Omega):=\left\{u\in L_{\Phi_i}(\Omega);\;\frac{\partial u}
{\partial x_i}\in L_{\Phi_i}(\Omega),\;i=1,...,N\right\}\,.$$
These are Banach spaces with respect to the norms
$$\|u\|_{1,\Phi_i}:=\|u\|_{\Phi_i}+\||\nabla u|\|_{\Phi_i},\;i=1,2\,.$$
We also define the Orlicz-Sobolev spaces $W_0^1L_{\Phi_i}(\Omega)$,
$i=1,2$, as the closure of $C_0^\infty(\Omega)$ in
$W^1L_{\Phi_i}(\Omega)$. By Lemma 5.7 in \cite{G} we obtain that on
$W_0^1L_{\Phi_i}(\Omega)$, $i=1,2$, we may consider some equivalent
norms
$$\|u\|_i:=\||\nabla u|\|_{\Phi_i}.$$
For an easier manipulation of the spaces defined above, we define
$$(\phi_i)_0:=\inf_{t>0}\frac{t\phi_i(t)}{\Phi_i(t)}\;\;{\rm and}\;\;(\phi_i)^0:=\sup_{t>0}\frac{t\phi_i(t)}{\Phi_i(t)},
\;i\in\{1,2\}\,.$$ 
In this paper we assume that for each
$i\in\{1,2\}$ we have
\begin{equation}\label{acdc0}1<(\phi_i)_0\leq\frac{t\phi_i(t)}{\Phi_i(t)}\leq(\phi_i)^0<\infty,\;\;\;\forall\;t\geq
0\,.\end{equation} The above relation implies that each $\Phi_i$,
$i\in\{1,2\}$, satisfies the $\Delta_2$-condition, i.e.
\begin{equation}\label{acdc}
\Phi_i(2t)\leq K\Phi_i(t),\;\;\;\forall\;t\geq 0\,,
\end{equation}
where $K$ is a positive constant (see  \cite[Proposition
2.3]{fourier}).

On the other hand, the following relations hold true
\begin{equation}\label{zz}
\|u\|_i^{(\phi_i)^0}\leq\int_\Omega\Phi_i(|\nabla u|)\;dx\leq\|u\|_i^{(\phi_i)_0},\;\;\;\forall\;u\in W_0^1L_{\Phi_i}(\Omega)\;
{\rm with}\;\|u\|_i<1,\;\;\;i=1,2\,,
\end{equation}
\begin{equation}\label{zz1}
\|u\|_i^{(\phi_i)_0}\leq\int_\Omega\Phi_i(|\nabla u|)\;dx\leq\|u\|_i^{(\phi_i)^0},\;\;\;\forall\;u\in W_0^1L_{\Phi_i}(\Omega)\;
{\rm with}\;\|u\|_i>1,\;\;\;i=1,2\,,
\end{equation}
(see, e.g. \cite[Lemma 1]{AA}).

Furthermore, in this paper we assume that for each
$i\in\{1,2\}$ the function $\Phi_i$ satisfies the following
condition
\begin{equation}\label{acdc1}
{\rm the}\;{\rm function}\;[0,\infty)\ni
t\rightarrow\Phi_i(\sqrt{t})\;{\rm is}\;{\rm convex}\,.
\end{equation}
Conditions \eq{acdc} and \eq{acdc1} assure that for each
$i\in\{1,2\}$ the Orlicz spaces $L_{\Phi_i}(\Omega)$ are
uniformly convex spaces and thus, reflexive Banach spaces (see
\cite[Proposition 2.2]{fourier}). That fact implies that also the
Orlicz-Sobolev spaces $W_0^1L_{\Phi_i}(\Omega)$,
$i\in\{1,2\}$, are reflexive Banach spaces.
\smallskip 

\noindent{\bf Remark 1.} We point out certain examples of
functions $\phi:\RR\rightarrow\RR$ which are odd, increasing
homeomorphisms from $\RR$ onto $\RR$ and satisfy conditions
\eq{acdc0} and \eq{acdc1}. For more details the reader can consult
\cite[Examples 1-3, p. 243]{Clem2}.

1) Let
$$\phi(t)=p|t|^{p-2}t,\;\;\;\forall\;t\in\RR\,,$$
with $p>1$. For this function it can be proved that
$$(\phi)_0=(\phi)^0=p\,.$$
Furthermore, in this particular case the corresponding Orlicz space $L_\Phi(\Omega)$ is the classical Lebesgue space 
$L^p(\Omega)$ while the Orlicz-Sobolev space $W_0^1L_\Phi(\Omega)$ is the classical Sobolev space $W_0^{1,p}(\Omega)$. 
We will use the classical notations to denote the Orlicz-Sobolev spaces in this particular case. 

2) Consider
$$\phi(t)=\log(1+|t|^s)|t|^{p-2}t,\;\;\;\forall\;t\in\RR\,,$$
with $p$, $s>1$. In this case it can be proved that
$$(\phi)_0=p,\;\;\;\;\;(\phi)^0=p+s\,.$$

3) Let
$$\phi(t)=\frac{|t|^{p-2}t}{\log(1+|t|)},\;\;\;{\rm if}\;t\neq 0,\;\;\;\phi(0)=0\,,$$
with $p>2$. In this case we have
$$(\phi)_0=p-1,\;\;\;\;\;(\phi)^0=p\,.$$
\smallskip

Next, we recall some background facts concerning the variable exponent Lebesgue spaces. For more details we refer to the book by Musielak
\cite{M} and the papers by Edmunds et al. \cite{edm, edm2, edm3}, Kovacik and  R\'akosn\'{\i}k \cite{KR},
Mih\u ailescu and R\u adulescu \cite{RoyalSoc}, and Samko and Vakulov \cite{samko}.

Set
$$C_+(\overline\Omega)=\{h;\;h\in C(\overline\Omega),\;h(x)>1\;{\rm
for}\;
{\rm all}\;x\in\overline\Omega\}.$$
For any $h\in C_+(\overline\Omega)$ we define
$$h^+=\sup_{x\in\Omega}h(x)\qquad\mbox{and}\qquad h^-=
\inf_{x\in\Omega}h(x).$$
For any $q(x)\in C_+(\overline\Omega)$ we define the variable exponent
Lebesgue space $L^{q(x)}(\Omega)$ (see \cite{KR}).
On $L^{q(x)}(\Omega)$ we define the Luxemburg norm by the formula
$$|u|_{q(x)}=\inf\left\{\mu>0;\;\int_\Omega\left|
\frac{u(x)}{\mu}\right|^{q(x)}\;dx\leq 1\right\}.$$ We remember that
the variable exponent Lebesgue spaces are separable and reflexive
Banach spaces. If $0 <|\Omega|<\infty$ and $q_1$, $q_2$ are variable
exponents so that $q_1(x) \leq q_2(x)$ almost everywhere in $\Omega$
then there exists the continuous embedding
$L^{q_2(x)}(\Omega)\hookrightarrow L^{q_1(x)}(\Omega)$.

Let $L^{p^{'}(x)}(\Omega)$ denote the conjugate space of
$L^{p(x)}(\Omega)$, where $1/p(x)+1/p^{'}(x)=1$. For any $u\in
L^{p(x)}(\Omega)$ and $v\in L^{p^{'}(x)}(\Omega)$ the H\"older
type inequality
\begin{equation}\label{Hol}
\left|\int_\Omega uv\;dx\right|\leq\left(\frac{1}{p^-}+
\frac{1}{{p^{'}}^-}\right)|u|_{p(x)}|v|_{p^{'}(x)}
\end{equation}
holds true.

If $(u_n)$, $u\in L^{q(x)}(\Omega)$ then the following relations
hold true
\begin{equation}\label{trei}
|u|_{q(x)}>1\;\;\;\Rightarrow\;\;\;|u|_{q(x)}^{q^-}\leq\int_\Omega|u|^{q(x)}\;dx
\leq|u|_{q(x)}^{q^+}
\end{equation}
\begin{equation}\label{patru}
|u|_{q(x)}<1\;\;\;\Rightarrow\;\;\;|u|_{q(x)}^{q^+}\leq
\int_\Omega|u|^{q(x)}\;dx\leq|u|_{q(x)}^{q^-}
\end{equation}
\begin{equation}\label{patrubiss}
|u_n-u|_{q(x)}\rightarrow 0\;\;\;\Leftrightarrow\;\;\;\int_\Omega|u_n-u|^{q(x)}\;dx
\rightarrow 0.
\end{equation}
\bigskip

Now we can turn back to problem \eq{1}. We will study  problem \eq{1} when  $q_1$, $q_2$, $m:\overline\Omega\rightarrow(1,\infty)$ are
continuous functions  satisfying the following assumptions:
\begin{equation}\label{2}
1<(\phi_2)_0\leq(\phi_2)^0<q_2^-\leq q_2^+\leq m^-\leq m^+\leq q_1^-\leq q_1^+<(\phi_1)_0\leq(\phi_1)^0<N\,,
\end{equation}
\begin{equation}\label{3}
q_1^+< [(\phi_2)_0]^\star:=\frac{N(\phi_2)_0}{N-(\phi_2)_0},\;\;\;\forall\;x\in\overline\Omega\, ,
\end{equation}
 and the potential $V:\Omega\rightarrow\RR$ satisfies
\begin{equation}\label{4}
V\in L^{r(x)}(\Omega),\;\;{\rm with}\;r(x)\in C(\overline\Omega)\;{\rm and}\;r(x)>\frac{N}{m^-}\;
\forall\;x\in\overline\Omega\,.
\end{equation}
Condition \eq{2} which describes  the competition between the growth rates  involved in equation \eq{1}, actually, assures a
balance between them and thus, it represents the {\it key} of the present study. Such a balance
is essential since we are working on a non-homogeneous (eigenvalue) problem for which a minimization technique based on
the Lagrange Multiplier Theorem can not be applied in order to find (principal) eigenvalues (unlike the case offered by the
homogeneous operators). Thus, in the case of nonlinear non-homogeneous eigenvalue problems the classical theory used in the
homogeneous case does not work entirely, but some of its ideas can still be useful and some particular results can still be
obtained in some aspects while in other aspects entirely new phenomena can occur. To focus on our case, condition \eq{2}
together with conditions \eq{3} and \eq{4} imply
$$\lim_{\|u\|_{1}\rightarrow 0}\di\frac{\di\int_\Omega\Phi_1(|\nabla u|)\;dx+\di\int_\Omega\Phi_2(|\nabla u|)\;dx+
\di\int_\Omega\frac{V(x)}{m(x)}|u|^{m(x)}\;dx}{\di\int_\Omega\frac{1}{q_1(x)}|u|^{q_1(x)}\;dx+
\di\int_\Omega\frac{1}{q_2(x)}|u|^{q_2(x)}\;dx}=\infty$$ and
$$\lim_{\|u\|_{1}\rightarrow\infty}\di\frac{\di\int_\Omega\Phi_1(|\nabla u|)\;dx+\di\int_\Omega\Phi_2(|\nabla u|)\;dx+
\di\int_\Omega\frac{V(x)}{m(x)}|u|^{m(x)}\;dx}{\di\int_\Omega\frac{1}{q_1(x)}|u|^{q_1(x)}\;dx+
\di\int_\Omega\frac{1}{q_2(x)}|u|^{q_2(x)}\;dx}=\infty\,.$$  In other words, the
absence of homogeneity is balanced by the behavior (actually, the
blow-up) of the Rayleigh quotient associated to problem \eq{1} in
the origin and at infinity. The consequences of the above remarks
is that the infimum of the Rayleigh quotient associated to problem
\eq{1} is a real number, i.e.
\begin{equation}\label{steamare}
\inf\limits_{u\in W_0^1L_{\Phi_1}(\Omega)\setminus\{0\}}\di\frac{\di\int_\Omega\Phi_1(|\nabla u|)\;dx+\di\int_\Omega\Phi_2(|\nabla u|)\;dx
+\di\int_\Omega\frac{V(x)}{m(x)}|u|^{m(x)}\;dx}{\di\int_\Omega\frac{1}{q_1(x)}|u|^{q_1(x)}\;dx+
\di\int_\Omega\frac{1}{q_2(x)}|u|^{q_2(x)}\;dx}\in\RR\,,
\end{equation}
and it will be attained for a function  $u_0\in W_0^{1,p_1(x)}(\Omega)\setminus\{0\}$. Moreover, the
value in \eq{steamare} represents an eigenvalue of problem \eq{1} with the corresponding eigenfunction $u_0$.
However, at this stage we can not say if the eigenvalue described above is the lowest eigenvalue of problem \eq{1}
or not, even if we are able to show that any $\lambda$ small enough is not an eigenvalue of \eq{1}. For the moment
this rests an open question. On the other hand, we can prove that any $\lambda$ superior to the value given
by relation \eq{steamare} is also an eigenvalue of problem \eq{1}. Thus, we conclude that problem \eq{1} possesses a
continuous family of eigenvalues.

Related with the above ideas we will also discuss the {\it optimization} of the eigenvalues described by relation
\eq{steamare} with respect to the potential $V$, providing that $V$ belongs to a bounded, closed and convex
subset of $L^{r(x)}(\Omega)$ (where $r(x)$ is given by relation \eq{4}). By optimization we understand
the existence of some potentials $V_\star$ and $V^\star$ such that the eigenvalue described in relation \eq{steamare}
is minimal or maximal with respect to the set where $V$ lies. The results that we will obtain in the context of optimization
of eigenvalues are motivated by the above advances in this field in the case of homogeneous (linear or nonlinear) eigenvalue
problems. We refer mainly to the studies in Asbaugh-Harrell \cite{AH}, Egnell \cite{egn} and Bonder-Del Pezzo \cite{bonder} where
different optimization problems of the principal eigenvalue of some homogeneous operators were studied.

\section{The main results}

By relation \eq{2} it follows that
$W_0^{1}L_{\Phi_1}(\Omega)$ is continuously embedded in
$W_0^{1}L_{\Phi_2}(\Omega)$ (see, e.g. \cite[Lemma 2]{AA}). Thus,  problem \eq{1} will be analyzed in the space 
$W_0^{1}L_{\Phi_1}(\Omega)$.

We say that $\lambda\in\RR$ is an {\it eigenvalue} of problem \eq{1} if
there exists $u\in W_0^{1}L_{\Phi_1}(\Omega)
\setminus\{0\}$ such that
$$\int_{\Omega}(a_1(|\nabla u|)+a_2(|\nabla u|))\nabla u\nabla v\;dx+
\int_\Omega V(x)|u|^{m(x)-2}uv\;dx-\lambda\int_{\Omega}
(|u|^{q_1(x)-2}+|u|^{q_2(x)-2})uv\;dx=0\,,$$ for all $v\in W_0^{1}L_{\Phi_1}(\Omega)$.
We point out that if $\lambda$ is an eigenvalue of problem \eq{1}
then the corresponding {\it eigenfunction} $u\in
W_0^{1}L_{\Phi_1}(\Omega)\setminus\{0\}$ is a {\it weak solution} of
problem \eq{1}.

For each potential $V\in L^{r(x)}(\Omega)$ we define
$$A(V):=\inf\limits_{u\in W_0^{1}L_{\Phi_1}(\Omega)\setminus\{0\}}
\di\frac{\di\int_\Omega\Phi_1(|\nabla u|)\;dx+\di\int_\Omega\Phi_2(|\nabla u|)\;dx
+\di\int_\Omega\frac{V(x)}{m(x)}|u|^{m(x)}\;dx}{\di\int_\Omega\frac{1}{q_1(x)}|u|^{q_1(x)}\;dx+
\di\int_\Omega\frac{1}{q_2(x)}|u|^{q_2(x)}\;dx}\,$$
and
$$B(V):=\inf\limits_{u\in W_0^{1}L_{\Phi_1}(\Omega)\setminus\{0\}}
\di\frac{\di\int_\Omega a_1(|\nabla u|)|\nabla u|^2\;dx+\di\int_\Omega 
a_2(|\nabla u|)|\nabla u|^2\;dx\;dx+\di\int_\Omega{V(x)}|u|^{m(x)}\;dx}{\di\int_\Omega|u|^{q_1(x)}\;dx+
\di\int_\Omega|u|^{q_2(x)}\;dx}\,.$$
Thus, we can define a functions $A,\; B:L^{r(x)}(\Omega)\rightarrow\RR$.

The first result of this paper is given by the following theorem.
\begin{teo}\label{t1}
Assume that conditions \eq{2}, \eq{3} and \eq{4} are fulfilled. Then $A(V)$ is an eigenvalue of problem \eq{1}.
Moreover, there exists $u_V \in W_0^{1}L_{\Phi_1}(\Omega)\setminus\{0\}$ an eigenfunction corresponding to
the eigenvalue $A(V)$ such that
$$A(V)=\di\frac{\di\int_\Omega\Phi_1(|\nabla u_V|)\;dx+\di\int_\Omega\Phi_2(|\nabla u_V|)\;dx
+\di\int_\Omega\frac{V(x)}{m(x)}|u_V|^{m(x)}\;dx}{\di\int_\Omega\frac{1}{q_1(x)}|u_V|^{q_1(x)}\;dx+
\di\int_\Omega\frac{1}{q_2(x)}|u_V|^{q_2(x)}\;dx}\,.$$
Furthermore, $B(V)\leq A(V)$, each $\lambda\in(A(V),\infty)$ is an eigenvalue of problem \eq{1}, while each
$\lambda\in(-\infty, B(V))$ is not an eigenvalue of problem \eq{1}.
\end{teo}
Next, we will show that on each convex, bounded and closed subset of $L^{r(x)}(\Omega)$ the function $A$
defined above is bounded from below and attains its minimum. The result is the following:
\begin{teo}\label{t2}
Assume that conditions \eq{2}, \eq{3} and \eq{4} are fulfilled. Assume that $S$ is a convex, bounded and closed
subset of $L^{r(x)}(\Omega)$. Then there exists $V_\star\in S$ which minimizes $A(V)$ on $S$, i.e.
$$A(V_\star)=\inf_{V\in S}A(V)\,.$$
\end{teo}
Finally, we will focus our attention on the particular case when the set $S$ from Theorem \ref{t2} is a
ball in $L^{r(x)}(\Omega)$. Thus, we will denote each closed ball centered in the origin of radius $R$ from $L^{r(x)}(\Omega)$
by $\overline B_R(0)$, i.e.
$$\overline B_R(0):=\{u\in L^{r(x)}(\Omega);\;\;|u|_{r(x)}\leq R\}\,.$$
By Theorem \ref{t2} we can define the function $A_\star:[0,\infty)\rightarrow\RR$ by
$$A_\star(R)=\min_{V\in \overline B_R(0)}A(V)\,.$$
Our result on the function $A_\star$ is given by the following theorem:
\begin{teo}\label{t3}
a) The function $A_\star$ is not constant and  decreases monotonically.

\noindent b)  The function $A_\star$ is continuous.
\end{teo}
On the other hand, we point out that  similar results as those of Theorems \ref{t2} and \ref{t3} can be obtained
if we notice that on each convex, bounded and closed subset of $L^{r(x)}(\Omega)$ the function $A$
defined in Theorem \ref{t1} is also bounded from above and attains its maximum. It is also easy to
remark that we can define a function $A^\star:[0,\infty)\rightarrow\RR$ by
$$A^\star(R)=\max_{V\in \overline B_R(0)}A(V)\,,$$
which has similar properties as $A_\star$.

\section{Proof of Theorem \ref{t1}}
Let $X$ denote the generalized Sobolev space
$W_0^{1}L_{\Phi_1}(\Omega)$. Relation \eq{2} and similar arguments as those used in \cite[Lemma 1]{AA} combined 
with \cite[Lemma 8.12(b)]{A} and with the Rellich-Kondrachov theorem we deduce that
\begin{equation}\label{embedding}
W_0^1L_{\Phi_1}(\Omega)\subset W_0^1L_{\Phi_2}(\Omega)\subset W_0^{1,(\phi_2)_0}(\Omega)\hookrightarrow L^{q_1^+}(\Omega)\subset
L^{q_1(x)}(\Omega)\subset L^{m(x)}(\Omega)\subset L^{q_2(x)}(\Omega)\,,
\end{equation}
where we denoted by $\subset$ a {\it continuous} embedding while by $\hookrightarrow$ we denoted a {\it compact} embedding.

Define the functionals $J_V$,  $I:X\rightarrow\RR$ by
$$J_V(u)=\int_{\Omega}\Phi_1(|\nabla u|)\;dx+\int_{\Omega}\Phi_2(|\nabla u|)\;dx
+\int_\Omega\frac{V(x)}{m(x)}|u|^{m(x)}\;dx\,,$$
$$I(u)=\int_\Omega\frac{1}{q_1(x)}|u|^{q_1(x)}\;dx+\int_\Omega\frac{1}{q_2(x)}|u|^{q_2(x)}\;dx\,.$$
Relation \eq{embedding} assures that the functionals defined above are well-defined. We notice that for any 
$V$ satisfying condition \eq{4} we have
$$J_V(u)=J_0(u)+\int_\Omega\frac{V(x)}{m(x)}|u|^{m(x)}\;dx,\;\;\;\forall\;u\in X\,,$$
where $J_0$ is obtained in the case when $V=0$ in $\Omega$.

Standard arguments imply that $J_V,\;I\in C^1(X,\RR)$ and for all
$u,\;v\in X$,
$$\langle J_V^{'}(u),v\rangle=\int_{\Omega}(a_1(|\nabla u|)+a_2(|\nabla u|))
\nabla u\nabla v\;dx+\int_\Omega V(x)|u|^{m(x)-2}uv\;dx\,,$$
$$\langle I^{'}(u),v\rangle=\int_\Omega|u|^{q_1(x)-2}uv\;dx+\int_\Omega|u|^{q_2(x)-2}uv\;dx\,.$$

In order to prove Theorem \ref{t1} we first establish some auxiliary results.
\begin{lemma}\label{l1}
Assume conditions \eq{2}, \eq{3} and \eq{4} are fulfilled. Then for each $\epsilon>0$ there exists
$C_\epsilon>0$ such that
$$
\left|\int_\Omega\frac{V(x)}{m(x)}|u|^{m(x)}\;dx\right|\leq\epsilon
\int_{\Omega}\left(\Phi_1(|\nabla u|)+\Phi_2(|\nabla u|)\right)\;dx+
C_\epsilon|V|_{r(x)}\int_\Omega(|u|^{m^-}+|u|^{m^+})\;dx\,,
$$
for all $u\in X$.
\end{lemma}
\proof
First, we point out that since $r(x)>r^{-}$ on $\overline\Omega$ it follows that $L^{r(x)}(\Omega)\subset L^{r^-}(\Omega)$.
On the other hand,  since $r(x)>\frac{N}{m^-}$ for each $x\in\overline\Omega$ it follows that $r^->\frac{N}{m^-}$. Thus,
we infer that $V\in L^{r^-}(\Omega)$ and  $r^->\frac{N}{m^-}$.

Now, let $\epsilon>0$ be fixed. We claim that there exists $D_\epsilon>0$ such that
\begin{equation}\label{e1}
\int_\Omega |V(x)|\cdot|u|^{m^-}\;dx\leq\epsilon\int_\Omega|\nabla u|^{m^-}\;dx+D_\epsilon|V|_{r^-}
\int_\Omega|u|^{m^-}\;dx,\;\;\;\forall\;u\in W^{1,m^-}_0(\Omega)\,.
\end{equation}
In order to establish \eq{e1} we show first that for each $s\in(1,\frac{Nm^-}{N-m^-})$  there exists $D_\epsilon^{'}>0$ such that
\begin{equation}\label{1stea}
|v|_s\leq\epsilon|\;|\nabla v|\;|_{m^-}+D_\epsilon^{'}|v|_{m^-},\;\;\;\forall\;u\in W^{1,m^-}_0(\Omega)\,.
\end{equation}
Indeed, assuming by contradiction that relation \eq{1stea} does not hold true for each $\epsilon>0$. Then there exists
$\epsilon_0>0$ and a sequence $(v_n)\subset W^{1,m^-}_0(\Omega)$ such that $|v_n|_{s}=1$ and
$$\epsilon_0|\;|\nabla v_n|\;|_{m^-}+n|v_n|_{m^-}<1,\;\;\;\forall\;n\,.$$
Then it is clear that $(v_n)$ is bounded in $W^{1,m^-}_0(\Omega)$ and $|v_n|_{m^-}\rightarrow 0$. Thus, we deduce that
passing eventually to a subsequence we can assume that $v_n$ converges weakly to a function $v$  in $W^{1,m^-}_0(\Omega)$ and
actually $v=0$. Since
$s\in(1,\frac{Nm^-}{N-m^-})$ it follows by the Rellich-Kondrachov theorem that $W^{1,m^-}_0(\Omega)$ is compactly
embedded in $L^{s}(\Omega)$ and thus $v_n$ converges to $0$ in $L^s(\Omega)$. On the other hand, since
$|v_n|_{s}=1$ for each $n$ we deduce that $|v|_{s}=1$ and that is a contradiction.   We obtained that relation \eq{1stea} holds true.

Next, we point out that since $r^->\frac{N}{m^-}$ then $m^-\cdot{r^-}^{'}<\frac{Nm^-}{N-m^-}$, where ${r^-}^{'}=\frac{r^-}{r^-1}$.
Thus, by H\"older's inequality we have
$$\int_\Omega |V(x)|\cdot|u|^{m^-}\;dx\leq|V|_{r^-}\cdot|u|_{m^-\cdot{r^-}^{'}}^{m^-},\;\;\;\forall\;u\in W^{1,m^-}_0(\Omega)\,.$$
Combining the last inequality with relation \eq{1stea} we infer that relation \eq{e1} holds true.

Similar arguments as those used in the proof of relation \eq{e1} combined with the fact that since $r^->\frac{N}{m^-}$
 we also have $r^->\frac{N}{m^+}$ imply that there exists $D_\epsilon^{''}$
\begin{equation}\label{e2}
\int_\Omega |V(x)|\cdot|u|^{m^+}\;dx\leq\epsilon\int_\Omega|\nabla u|^{m^+}\;dx+D_\epsilon^{''}|V|_{r^-}
\int_\Omega|u|^{m^+}\;dx,\;\;\;\forall\;u\in W^{1,m^+}_0(\Omega)\,.
\end{equation}
Using relation \eq{2} we deduce that $m^-\leq m^+<(\phi_1)_0$  and thus,
implies that $W_0^{1,(\phi_1)_0}(\Omega)\subset W_0^{1,m^{\pm}}(\Omega)$. On the other hand, similar arguments as those 
used in the proof of \cite[Lemma 2]{AA} show that $W_0^1L_{\Phi_1}(\Omega)\subset W_0^{1,(\phi_1)_0}(\Omega)$. 
The above facts imply that relations \eq{e1} and \eq{e2} hold true for any $u\in X$. Moreover,
in the right hand sides of inequalities \eq{e1} and \eq{e2} we can take $|V|_{r(x)}$ instead of $|V|_{r^-}$ since
$L^{r(x)}(\Omega)$ is continuously embedded in $L^{r^-}(\Omega)$ via inequality \eq{Hol}.

Finally, we point out that since by \eq{2} we have $(\phi_2)^0<m^-\leq m(x)\leq m^+<(\phi_1)_0$ for each $x\in\overline\Omega$ we deduce
that
\begin{equation}\label{e3}
\left|\int_\Omega\frac{V(x)}{m(x)}|u|^{m(x)}\;dx\right|\leq\frac{1}{m^-}\int_\Omega|V(x)|\cdot(|u|^{m^-}+|u|^{m^+})\;dx,
\;\;\;\forall\;u\in X
\end{equation}
and
\begin{equation}\label{e4}
\int_\Omega(|\nabla u|^{m^-}+|\nabla u|^{m^+})\;dx\leq 
\int_{\Omega}\left(|\nabla u|^{(\phi_2)^0}\;dx+|\nabla u|^{(\phi_1)_0}\right)\;dx,\;\;\;\forall\;u\in X\,.
\end{equation}
Relations \eq{e1}, \eq{e2}, \eq{e3}, \eq{e4}, \eq{2} and \eq{embedding} and \cite[Lemma 3]{AA} 
lead to the idea that Lemma \ref{l1} holds true.  \endproof

\begin{lemma}\label{l2}
The following relations hold true:
\begin{equation}\label{10}
\lim_{\|u\|_1\rightarrow\infty}\frac{J_V(u)}{I(u)}=\infty
\end{equation}
and
\begin{equation}\label{11}
\lim_{\|u\|_1\rightarrow 0}\frac{J_V(u)}{I(u)}=\infty.
\end{equation}
\end{lemma}
\proof
First, we point out that by \eq{2} $q_2(x)<m^{\pm}< q_1(x)$ for any $x\in\overline\Omega$. Thus, it is clear that
$$|u(x)|^{m^-}+|u(x)|^{m^+}\leq 2(|u(x)|^{q_1(x)}+|u(x)|^{q_2(x)}),\;\;\;\forall\;x\in\overline\Omega\;{\rm and}\;\forall\;
u\in X\,.$$
Integrating over $\Omega$ the above inequality we infer that
\begin{equation}\label{163276}
\frac{\di\int_\Omega(|u|^{m^-}+|u|^{m^+})\;dx}{\di\int_\Omega(|u|^{q_1(x)}+|u|^{q_2(x)})\;dx}\leq 2,\;\;\;\forall\;u\in X\,.
\end{equation}
Using Lemma \ref{l1}  we find that for an $\epsilon\in(0,1)$ there exists $C_\epsilon>0$ such that
$$\frac{J_V(u)}{I(u)}\geq\frac{{(1-\epsilon)}\di\int_\Omega(\Phi_1(|\nabla u|)+\Phi_2(|\nabla u|))\;dx-C_\epsilon
|V|_{r(x)}\di\int_\Omega(|u|^{m^-}+|u|^{m^+})\;dx}{\di\frac{1}{q_2^-}\di\int_\Omega(|u|^{q_1(x)}+|u|^{q_2(x)})\;dx}\,,$$
for any $u\in X$.

By the above inequality and relation \eq{163276} we deduce that there exist some positive constants $\beta>0$ and $\gamma>0$
such that
\begin{equation}\label{pietrosu}
\frac{J_V(u)}{I(u)}\geq\frac{\beta\di\int_\Omega(\Phi_1(|\nabla u|)+\Phi_2(|\nabla u|))\;dx}
{\di\int_\Omega(|u|^{q_1(x)}+|u|^{q_2(x)})\;dx}-\gamma|V|_{r(x)},\;\;\;\forall\;u\in X\,.
\end{equation}
For any $u\in X$ with $\|u\|_1>1$  relation \eq{pietrosu}  implies
$$\frac{J_V(u)}{I(u)}\geq\frac{\beta\di\int_\Omega\Phi_1(|\nabla u|)\;dx}
{|u|_{q_1^-}^{q_1^-}+|u|_{q_1^+}^{q_1^+}+|u|_{q_2^-}^{q_2^-}+|u|_{q_2^+}^{q_2^+}}-\gamma|V|_{r(x)},\;\;\;\forall\;u\in X\;{\rm with}\;
\|u\|_1>1\,.$$
Now, taking into account the continuous embedding of $X$ in $L^{q_i^\pm}(\Omega)$ for $i=1,2$ (given by relations \eq{2} and \eq{embedding}) 
and the result of relation \eq{zz1} we deduce the existence of a positive constant $\delta>0$ such that
$$\frac{J_V(u)}{I(u)}\geq\frac{\delta\|u\|_1^{(\phi_1)_0}}
{\|u\|_1^{q_1^-}+\|u\|_1^{q_1^+}+\|u\|_1^{q_2^-}+\|u\|_1^{q_2^+}}-\gamma|V|_{r(x)},\;\;\;\forall\;u\in X\;{\rm with}\;
\|u\|_1>1\,.$$
Since $(\phi_1)_0>q_1^+\geq q_1^-\geq q_2^+\geq q_2^-$, passing to the limit as $\|u\|_1\rightarrow\infty$ in the above inequality we deduce that
relation \eq{10} holds true.

Relation \eq{embedding} shows that the space $W_0^{1}L_{\Phi_1}(\Omega)$ is
continuously embedded in $W_0^{1}L_{\Phi_2}(\Omega)$. Thus, if
$\|u\|_1\rightarrow 0$ then $\|u\|_2\rightarrow 0$.

The above remarks enable us to affirm that for any $u\in X$ with $\|u\|_1<1$ small enough we have $\|u\|_2<1$.

Using again relation \eq{embedding} we deduce that $W_0^{1}L_{\Phi_2}(\Omega)$ is continuously
embedded in $L^{q_i^\pm}(\Omega)$ with $i=1,2$.
 It follows that there exist four positive constants
$d_{i1}$ and $d_{i2}$ with $i=1,2$ such that
\begin{equation}\label{12prim}
\|u\|_2\geq d_{i1}\cdot|u|_{q_i^+},\;\;\;\forall\;u\in W_0^{1}L_{\Phi_2}(\Omega)\;{\rm and}\;i=1,2
\end{equation}
and
\begin{equation}\label{13prim}
\|u\|_2\geq d_{i2}\cdot|u|_{q_i^-},\;\;\;\forall\;u\in W_0^{1}L_{\Phi_2}(\Omega)\;{\rm and}\;i=1,2\,.
\end{equation}
Thus, for any $u\in X$ with $\|u\|_1<1$ small enough, relation \eq{pietrosu} implies
$$\frac{J_V(u)}{I(u)}\geq\frac{\beta\di\int_\Omega\Phi_2(|\nabla u|)\;dx}
{|u|_{q_1^-}^{q_1^-}+|u|_{q_1^+}^{q_1^+}+|u|_{q_2^-}^{q_2^-}+|u|_{q_2^+}^{q_2^+}}-\gamma|V|_{r(x)}\,.$$
Next, relations \eq{zz}, \eq{12prim}, \eq{13prim} yield that there exists  a constant $\xi>0$ such that
$$\frac{J_V(u)}{I(u)}\geq\frac{\xi\|u\|_2^{(\phi_2)^0}}
{\|u\|_2^{q_1^-}+\|u\|_2^{q_1^+}+\|u\|_2^{q_2^-}+\|u\|_2^{q_2^+}}-\gamma|V|_{r(x)}\,,$$
for any $u\in X$ with $\|u\|_1<1$ small enough.
Since $(\phi_2)^0<q_2^-\leq q_2^+\leq q_1^-\leq q_1^+$, passing to the limit as
$\|u\|_1\rightarrow 0$ (and thus, $\|u\|_2\rightarrow 0$) in the
above inequality we deduce that relation \eq{11} holds true. The
proof of Lemma \ref{l2} is complete.   \endproof

\noindent{\bf Remark 2.} We point out that by relation \eq{pietrosu} and using similar arguments as in the proof of
Theorem 1 (Step 1) in \cite{AA} we can find that for $V$ given and satisfying \eq{4} the quotient
$\frac{J_V(u)}{I(u)}$ is bounded from below for $u\in X\setminus\{0\}$, i.e. $A(V)$ is a real number. Similarly, 
it can be proved that $B(V)$ is also a real number.

\begin{lemma}\label{l3}
There exists $u\in X\setminus\{0\}$ such that $\frac{J_V(u)}{I(u)}=A(V)$.
\end{lemma}
\proof Let $(u_n)\subset X\setminus\{0\}$ be a minimizing
sequence for $A(V)$, that is,
\begin{equation}\label{stea}
\lim_{n\rightarrow\infty}\frac{J_V(u_n)}{I(u_n)}=A(V)\,.
\end{equation}
By relation \eq{10} it is clear that $\{u_n\}$ is bounded in $X$.
Since $X$ is reflexive it follows that there exists $u\in X$ such
that, up to a subsequence, $(u_n)$ converges weakly to $u$ in
$X$. On the other hand, similar arguments as those used in the
proof of \cite[Theorem 2]{jmaa} (see also \cite[Step 3]{AA}) show that the functional $J_0$ (obtained for $V=0$ on $\Omega$)
is weakly lower semi-continuous. Thus, we find
\begin{equation}\label{14}
\liminf_{n\rightarrow\infty}J_0(u_n)\geq J_0(u)\,.
\end{equation}
By the compact embedding theorem for Sobolev spaces and assumptions \eq{2}, \eq{3} and \eq{4} it follows that $X$ is
compactly embedded in $L^{\sigma(x)}(\Omega)$ (where $\sigma(x)=m(x)\cdot{r(x)}/({r(x)-1}$))
and $L^{q_i(x)}(\Omega)$ with $i=1,2$. Thus, $(u_n)$ converges strongly in $L^{\sigma(x)}(\Omega)$
and $L^{q_i(x)}(\Omega)$ with $i=1,2$. Then, by relations \eq{Hol} and
\eq{embedding} it follows that
\begin{equation}\label{15}
\lim_{n\rightarrow\infty}I(u_n)=I(u)
\end{equation}
and
\begin{equation}\label{15prim}
\lim_{n\rightarrow\infty}\int_\Omega V(x)|u_n|^{m(x)}\;dx=\int_\Omega V(x)|u|^{m(x)}\;dx\,.
\end{equation}
Relations \eq{14}, \eq{15} and \eq{15prim} imply that if $u\not\equiv 0$ then
$$\frac{J_V(u)}{I(u)}=A(V)\,.$$
Thus, in order to conclude that the lemma holds true it is enough
to show that $u$ is not trivial. Assume by contradiction the
contrary. Then $u_n$ converges weakly to $0$ in $X$ and strongly
in $L^{s(x)}(\Omega)$ for any $s(x)\in C(\overline\Omega)$ with $1<s(x)<\frac{N(\phi_1)_0}{N-(\phi_1)_0}$ on $\overline\Omega$.
In other words, we will have
\begin{equation}\label{victor}
\lim_{n\rightarrow\infty}I(u_n)=0\,,
\end{equation}
and
\begin{equation}\label{lucian}
\lim_{n\rightarrow\infty}\int_\Omega V(x)|u_n|^{m(x)}\;dx=0\,.
\end{equation}
Letting $\epsilon\in(0,|A(V)|)$ be fixed by relation \eq{stea} we deduce that for $n$ large enough we have
$$|J_V(u_n)-A(V)I(u_n)|<\epsilon I(u_n)\,,$$
or
$$(|A(V)|-\epsilon)I(u_n)<J_V(u_n)<(|A(V)|+\epsilon)I(u_n)\,.$$
Passing to the limit in the above inequalities and taking into account that relation \eq{victor} holds true we find
$$\lim_{n\rightarrow\infty}J_V(u_n)=0\,.$$
Next, by relation \eq{lucian} we get
$$\lim_{n\rightarrow\infty}J_0(u_n)=0\,.$$
That fact combined with relation \eq{zz} implies that actually $u_n$ converges strongly to $0$ in $X$, i.e.
$\lim_{n\rightarrow\infty}\|u_n\|_1=0$. By this information and relation \eq{11} we get
 $$\lim_{n\rightarrow\infty}\frac{J_V(u_n)}{I(u_n)}=\infty,$$
 and this is a contradiction. Thus, $u\not\equiv 0$.
 The proof of Lemma \ref{l3} is complete.  \endproof

 By Lemma \ref{l3} we conclude that there exists $u\in X\setminus\{0\}$ such that
\begin{equation}\label{16}
\frac{J_V(u)}{I(u)}=A(V)=\inf_{w\in X\setminus\{0\}}\frac{J_V(w)}{I(w)}\,.
\end{equation}
Then, for any $w\in X$ we have
$$\frac{d}{d\epsilon}\frac{J_V(u+\epsilon w)}{I(u+\epsilon w)}\left|_{\epsilon=0}=0\right.\,.$$
A simple computation yields
\begin{equation}\label{17}
\langle J_V^{'}(u),w\rangle I(u)-J_V(u)\langle I^{'}(u),w\rangle=0\,,
\end{equation}
for all $w\in X$.
Relation \eq{17} combined with the fact that $J_V(u)=A(V)\cdot I(u)$
and $I(u)\neq 0$ implies the fact that $A(V)$ is an eigenvalue
of problem \eq{1}.
\medskip

Next, we show that any
$\lambda\in(A(V),\infty)$ is an eigenvalue of problem \eq{1}.

Let $\lambda\in(A(V),\infty)$ be arbitrary but fixed. Define $T_{V,\lambda}:X\rightarrow\RR$ by
$$T_{V,\lambda}(u)=J_V(u)-\lambda I(u)\,.$$
Clearly, $T_{V,\lambda}\in C^1(X,\RR)$ with
$$\langle T_{V,\lambda}^{'}(u),v\rangle=\langle J_V^{'}(u),v\rangle-\lambda\langle I^{'}(u),
v\rangle,\;\;\;\forall\;u\in X.$$
Thus, $\lambda$ is an eigenvalue of problem \eq{1} if and only if there exists $u_\lambda\in X\setminus\{0\}$ a
critical point of $T_{V,\lambda}$.

With similar arguments as in the proof of relation \eq{10} we can
show that $T_{V,\lambda}$ is coercive, i.e.
$\lim_{\|u\|\rightarrow\infty}T_{V,\lambda}(u)=\infty$. On the other
hand, as we have already remarked, similar arguments as those used
in the proof of   \cite[Theorem 2]{jmaa} show that the
functional $T_{V,\lambda}$ is weakly lower semi-continuous. These two
facts enable us to apply  \cite[Theorem 1.2]{S} in order to prove
that there exists $u_\lambda\in X$ a global minimum point of
$T_{V,\lambda}$ and thus, a critical point of $T_{V,\lambda}$. It is enough to show that
$u_\lambda$ is not trivial. Indeed, since $A(V)=\inf_{u\in
X\setminus\{0\}}\frac{J_V(u)}{I(u)}$ and $\lambda>A(V)$ it
follows that there exists $v_\lambda\in X$ such that
$$J_V(v_\lambda)<\lambda I(v_\lambda)\,,$$
or
$$T_{V,\lambda}(v_\lambda)<0\,.$$
Thus,
$$\inf_{X}T_{V,\lambda}<0$$
and we conclude that $u_\lambda$ is a nontrivial critical point of
$T_{V,\lambda}$, or $\lambda$ is an eigenvalue of problem \eq{1}.

Finally, we prove that each $\lambda<B(V)$ is not an eigenvalue of problem \eq{1}. With
that end in view we assume by contradiction that there exists  $\lambda<B(V)$ an eigenvalue of problem \eq{1}.
It follows that there exists $u_\lambda\in X\setminus\{0\}$ such that
$$\langle J_V^{'}(u_\lambda),u_\lambda\rangle=\lambda\langle I^{'}(u_\lambda),u_\lambda\rangle\,.$$
Since $u_\lambda\neq 0$ we have $\langle I^{'}(u_\lambda),u_\lambda\rangle>0$.
Using that fact and the definition of $B(V)$ it follows that the  following relation holds true
$$\langle J_V^{'}(u_\lambda),u_\lambda\rangle=\lambda\langle I^{'}(u_\lambda),u_\lambda\rangle
<B(V)\langle I^{'}(u_\lambda),u_\lambda\rangle\leq\langle J_V^{'}(u_\lambda),u_\lambda\rangle \,.$$
Obviously, this is a contradiction. We deduce that each $\lambda\in(-\infty,B(V))$ is not an eigenvalue of problem \eq{1}.
Furthermore, it is clear that $A(V)\geq B(V)$.

The proof of Theorem \ref{t1} is complete.  \endproof
\smallskip

\noindent{\bf Remark 3.} We point out that in the case when $V=0$ in $\Omega$ the same arguments as in the
proof of Theorem 1 (Step 1) in \cite{AA} assure that $A(0)>0$.

\section{Proof of Theorem \ref{t2}}
Let $S$ be a convex, bounded and closed subset of $L^{r(x)}(\Omega)$ and
$$A_\star:=\inf_{V\in S}A(V)\,.$$
Clearly, relation \eq{pietrosu} assures that $A_\star$ is finite.

On the other hand, let $(V_n)\subset S$ be a minimizing sequence for $A_\star$, i.e.
$$A(V_n)\rightarrow A_\star,\;\;\;{\rm as}\;n\rightarrow\infty\,.$$
Obviously, $(V_n)$ is a bounded sequence and thus, there exists $V_\star\in L^{r(x)}(\Omega)$ such that
$V_n$ converges weakly to $V_\star$ in $L^{r(x)}(\Omega)$. Moreover, since $S$ is convex and closed it is
also weakly closed (see, e.g.,  Brezis \cite[Theorem III.7]{B}) and consequently $V_\star\in S$.

Next, we will show that $A(V_\star)=A_\star$.

Indeed, by Theorem \ref{t1} we deduce that for each positive integer $n$ there exists $u_n\in X\setminus\{0\}$ such that
\begin{equation}\label{dana}\frac{J_{V_n}(u_n)}{I(u_n)}=A(V_n)\,.\end{equation}
Since $(A(V_n))$ is a bounded sequence and by relation \eq{pietrosu} we have
$$\frac{J_{V_n}(u_n)}{I(u_n)}\geq\beta\frac{J_{0}(u_n)}{I(u_n)}-C,\;\;\;{\rm for}\;{\rm any}\;n\,,$$
where $C$ is a positive constant, we infer that $(u_n)$ is bounded in $X$ and it can not
contain a subsequence converging to $0$  (otherwise we obtain a contradiction by applying
Lemma \ref{l2}). Thus, there exists $u_0\in X\setminus\{0\}$ such that $(u_n)$ converges weakly to $u_0$ in $X$.  Using relation 
\eq{3} (and thus,  $W_0^1L_{\Phi_1}(\Omega)\subset W_0^{1,(\phi_1)_0}(\Omega)$) and the Rellich-Kondrachov theorem we deduce 
that $(u_n)$ converges strongly to $u_0$ in $L^{s(x)}(\Omega)$ for any $s(x)\in C(\overline\Omega)$ satisfying 
$1<s(x)<\frac{N(\phi_1)_0}{N-(\phi_1)_0}$ for any $x\in\overline\Omega$. In particular, using conditions \eq{2}, \eq{3} and \eq{4} we get  
that $(u_n)$ converges to $u_0$ in $L^{m(x)}(\Omega)$ and in $L^{m(x)\cdot r^{'}(x)}(\Omega)$ where $r^{'}(x)=\frac{r(x)}{r(x)-1}$.
Using that information, inequality \eq{Hol} and the fact that $V_\star\in L^{r(x)}(\Omega)$ and $(V_n)$ is bounded in $L^{r(x)}(\Omega)$ we find
\begin{equation}\label{1770815}
\lim_{n\rightarrow\infty}\int_\Omega\frac{V_\star(x)}{m(x)}|u_n|^{m(x)}\;dx=\int_\Omega\frac{V_\star(x)}{m(x)}|u_0|^{m(x)}\;dx
\end{equation}
and
\begin{equation}\label{2810123}
\lim_{n\rightarrow\infty}\int_\Omega\left(\frac{V_n(x)}{m(x)}|u_n|^{m(x)}-\frac{V_n(x)}{m(x)}|u_0|^{m(x)}\right)\;dx=0\,.
\end{equation}
On the other hand, since $(V_n)$ converges weakly to $V_\star$ in $L^{r(x)}(\Omega)$ and $u_0\in L^{m(x)\cdot r^{'}(x)}(\Omega)$,
where $r^{'}(x)=\frac{r(x)}{r(x)-1}$, we deduce
\begin{equation}\label{2490301}
\lim_{n\rightarrow\infty}\int_\Omega\frac{V_n(x)}{m(x)}|u_0|^{m(x)}\;dx=\int_\Omega\frac{V_\star(x)}{m(x)}|u_0|^{m(x)}\;dx\,.
\end{equation}
Combining the equality
\begin{eqnarray*}
\int_\Omega\frac{V_\star(x)}{m(x)}|u_n|^{m(x)}\;dx-\int_\Omega\frac{V_n(x)}{m(x)}|u_n|^{m(x)}\;dx
&=&\int_\Omega\frac{V_\star(x)}{m(x)}|u_n|^{m(x)}\;dx-\int_\Omega\frac{V_\star(x)}{m(x)}|u_0|^{m(x)}\;dx\\
&+&\int_\Omega\frac{V_\star(x)}{m(x)}|u_0|^{m(x)}\;dx-\int_\Omega\frac{V_n(x)}{m(x)}|u_0|^{m(x)}\;dx\\
&+&\int_\Omega\frac{V_n(x)}{m(x)}|u_0|^{m(x)}\;dx-\int_\Omega\frac{V_n(x)}{m(x)}|u_n|^{m(x)}\;dx\,,
\end{eqnarray*}
with relations \eq{1770815}, \eq{2810123} and \eq{2490301} we get
\begin{equation}\label{1500318}
\lim_{n\rightarrow\infty}\int_\Omega\left(\frac{V_\star(x)}{m(x)}|u_n|^{m(x)}-\frac{V_n(x)}{m(x)}|u_n|^{m(x)}\right)\;dx=0\,.
\end{equation}
Since
$$A(V_\star)=\inf_{u\in X\setminus\{0\}}\frac{J_{V_\star}(u)}{I(u)}\,,$$
it follows that
$$A(V_\star)\leq\frac{J_{V_\star}(u_n)}{I(u_n)}\,.$$
Combining the above inequality and equality \eq{dana} we obtain
$$A(V_\star)\leq\frac{J_{V_\star}(u_n)-J_{V_n}(u_n)}{I(u_n)}+A(V_n)\,.$$
Taking into account the result of relation \eq{1500318}, the fact that $I(u_n)$ is bounded and does not converge to
$0$ and $(A(V_n))$ converges to $A_\star$ then passing to the limit as $n\rightarrow\infty$ in the last inequality we infer that
$$A(V_\star)\leq A_\star\,.$$
But using the definition of $A_\star$ and the fact that $V_\star\in S$ we conclude that actually
$$A(V_\star)= A_\star\,.$$
The proof of Theorem \ref{t2} is complete. \endproof

\section{Proof of Theorem \ref{t3}}
a) First, we show that  function $A_\star$ is not constant. Indeed, by Remark 3 we point out that
$A_\star(0)=A(0)>0$. On the other hand, by  \cite[Theorem 1]{AA} it follows that
$$\lambda_m:=\inf_{u\in X\setminus\{0\}}\di\frac{\di\int_\Omega\Phi_1(|\nabla u|)\;dx+\di\int_\Omega\Phi_2(|\nabla u|)\;dx}
{\di\int_\Omega\frac{1}{m(x)}|u|^{m(x)}\;dx}>0\,.$$
Moreover,  \cite[Lemma 5]{AA} implies that there exists $u_m\in X\setminus\{0\}$ such that
$$\lambda_m=\frac{\di\int_\Omega\Phi_1(|\nabla u_m|)\;dx+\di\int_\Omega\Phi_2(|\nabla u_m|)\;dx}
{\di\int_\Omega\frac{1}{m(x)}|u_m|^{m(x)}\;dx}\,.$$
Thus, taking $V_m(x)=-\lambda_m$ for all $x\in\Omega$ it is clear that $V_m\in L^\infty(\Omega)\subset L^{r(x)}(\Omega)$ and
$$\frac{J_{V_m}(u_m)}{I(u_m)}=0\,.$$
It follows that
$$A(V_m)\leq 0\,,$$
and we find
$$A_\star(\lambda_m)\leq 0\,.$$
We conclude that $A_\star$ is not constant. Furthermore, we point out that a similar proof as those
presented above can show that function $A_\star$ takes also negative values. To support that idea
we just notice that by  \cite[Theorem 1, Step 3]{AA} for each $\lambda>\lambda_m$ there exits $u_\lambda\in X\setminus\{0\}$
such that taking $V_\lambda=-\lambda$ for all $x\in\Omega$ we have
$$\frac{J_{V_\lambda}(u_\lambda)}{I(u_\lambda)}<0\,.$$

Next, we point out  that $A_\star$ decreases monotonically. Indeed, if
we consider $0\leq R_1<R_2$ then it is clear that $\overline B_{R_1}(0)\subset \overline B_{R_2}(0)$. Then
the definition of  function $A_\star$ implies $A_\star(R_1)\geq A_\star(R_2)$.
\smallskip

b) Finally, we show that the function $A_\star$ is continuous. Let $R>0$ and $t\in(0,R)$ be fixed.
We will verify that $\lim_{t\searrow 0}A_\star(R+t)=\lim_{t\searrow 0}A_\star(R-t)=A_\star(R)$.

First, we prove that $\lim_{t\searrow 0}A_\star(R+t)=A_\star(R)$.
By Theorem \ref{t3} a) we have
$$A_\star(R)\geq A_\star(R+t)\,.$$
Moreover, by Theorem \ref{t2} it follows that there exists $V_{R+t}\in\overline B_{R+t}(0)$ (i.e. $|V_{R+t}|_{r(x)}\leq R+t$)
such that
$$A(V_{R+t})=A_\star(R+t)\,.$$
Taking now $V_{R,t}:=\frac{R}{R+t}V_{R+t}$ we have
$$|V_{R,t}|_{r(x)}=\frac{R}{R+t}|V_{R+t}|_{r(x)}\leq R\,,$$
or $V_{R,t}\in\overline B_R(0)$. Therefore, obviously, we have $A(V_{R,t})\geq A_\star(R)$.

On the other hand, by Theorem \ref{t1} there exists $u_t\in X\setminus\{0\}$ such that
$$A(V_{R+t})=\frac{J_{V_{R+t}}(u_t)}{I(u_t)}\,.$$
Combining the above pieces of information we find
\begin{eqnarray*}
A_\star(R+t)=A(V_{R+t})&=&\frac{J_{V_{R+t}}(u_t)}{I(u_t)}\\
&=&\frac{J_{\frac{R+t}{R}\cdot V_{R,t}}(u_t)}{I(u_t)}\\
&=&\frac{R+t}{R}\cdot\frac{J_{V_{R,t}}(u_t)}{I(u_t)}-
\frac{t}{R}\cdot\di\frac{J_0(u_t)}{I(u_t)}\\
&\geq& \frac{R+t}{R}\cdot A_\star(R)-\frac{t}{R}\cdot\di\frac{J_0(u_t)}{I(u_t)}\,.
\end{eqnarray*}
On the other hand, by relation \eq{pietrosu} we have that for each $t\in(0,R)$ it holds
\begin{eqnarray*}
A_\star(R)\geq A_\star(R+t)=A(V_{R+t})&=&\frac{J_{V_{R+t}}(u_t)}{I(u_t)}\\
&\geq&\beta_1\cdot\frac{J_0(u_t)}{I(u_t)}-\gamma\cdot|V_{R+t}|_{r(x)}\\
&=&\beta_1\cdot\frac{J_0(u_t)}{I(u_t)}-\gamma\cdot 2R\,,
\end{eqnarray*}
where $\beta_1>0$ and $\gamma>0$ are real constants.

Combining the last two inequalities we deduce that
$$A_\star(R)\geq A_\star(R+t)\geq\frac{R+t}{R}\cdot A_\star(R)-\frac{t}{R}\cdot\frac{A_\star(R)+\gamma\cdot 2R}{\beta_1}\,,$$
for each $t\in(0,R)$.

We conclude that
$$\lim_{t\searrow 0}A_\star(R+t)=A_\star(R)\,.$$
In the following we argue that $\lim_{t\searrow 0}A_\star(R-t)=A_\star(R)$.

Obviously,
$$A_\star(R)\leq A_\star(R-t),\;\;\;\forall\;t\in(0,R)\,.$$
By Theorem \ref{t2} there exists $V_R\in\overline B_R(0)$ such that
$$A_\star(R)=A(V_R)\,.$$
Moreover, by Theorem \ref{t1} there exists $u_0\in X\setminus\{0\}$ such that
$$A(V_R)=\frac{J_{V_R}(u_0)}{I(u_0)}\,.$$
Define now
$$V_t:=\frac{R-t}{R} V_R,\;\;\;\forall\;t\in(0,R)\,.$$
Clearly, $V_t\in \overline B_{R-t}(0)$. Thus, it is clear that
$$\frac{J_{V_t}(u_0)}{I(u_0)}\geq A_\star(R-t),\;\;\;\forall\;t\in(0,R)\,.$$
Taking into account the above information we find
\begin{eqnarray*}
A_\star(R)=A(V_R)=\frac{J_{V_R}(u_0)}{I(u_0)}&=&
\frac{J_{\frac{R}{R-t}V_t}(u_0)}{I(u_0)}\\
&=&\frac{J_{V_t}(u_0)}{I(u_0)}+\frac{t}{R-t}\cdot\frac{\di\int_\Omega\frac{V_t(x)}{m(x)}|u_0|^{m(x)}\;dx}{I(u_0)}\\
&\geq&A_\star(R-t)+\frac{t}{R}\cdot\frac{\di\int_\Omega\frac{V_R(x)}{m(x)}|u_0|^{m(x)}\;dx}
{I(u_0)},\;\;\;\forall\;t\in(0,R)\,.
\end{eqnarray*}
We infer
$$\lim_{t\searrow 0}A_\star(R-t)=A_\star(R)\,.$$
It follows that  function  $A_\star$ is continuous. The proof of Theorem \ref{t3} is complete. \endproof
\medskip

\noindent{\bf Remark 4.} By Theorem \ref{t3} a) we get that $A_\star$ decreases monotonically. We notice that in
the particular case when $q_1(x)=m(x)=q_2(x)=q$ for each $x\in\overline\Omega$, where $q>1$ is a real number
for which conditions \eq{2}, \eq{3} and \eq{4} are fulfilled, the above quoted result can be improved, in the sense
that we can show that, actually, function $A_\star$ is strictly decreasing on $[0,\infty)$. Indeed, letting
$0\leq R_1<R_2$ be given, by Theorem \ref{t2} we deduce that there exists $V_1\in \overline B_{R_1}(0)$ such
that
$$A(V_1)=A_\star(R_1)\,.$$
Then for each real number $t\in(0,R_2-R_1)$ we have $V_1-t\in \overline B_{R_2}(0)$ since
$|V_1-t|_{r(x)}\leq|V_1|_{r(x)}+t\leq R_2$. Next, by Theorem \ref{t1} there exists $u_1\in X\setminus\{0\}$ such
that
$$A(V_1)=\frac{J_{V_1}(u_1)}{I(u_1)}\,.$$
Taking into account all the above remarks we infer
$$A_\star(R_1)-\frac{t}{2}=A(V_1)-\frac{t}{2}=\frac{J_{V_1}(u_1)}{I(u_1)}-\frac{t}{2}=\frac{J_{V_1-t}(u_1)}{I(u_1)}
\geq A(V_1-t)\geq A_\star(R_2)\,,$$
or
$$A_\star(R_1)>A_\star(R_2)\,.$$
In the end of this remark we consider that it is important to highlight the idea that the above proof supports the
fact that in the case when we manipulate homogeneous quantities we obtain better results than in the case
when we deal with non-homogeneous quantities.
\medskip

\noindent{\bf Remark 5.} We point out that by Theorem \ref{t3} b) we deduce that
$$A_\star(R)=\inf_{s\leq R}A_\star(s)\;\;\;{\rm and}\;\;\;A_\star(R)=\sup_{s\geq R}A_\star(s)\,.$$
\medskip

\noindent{\bf Remark 6.} We also point out that  function $A_\star$ can be used in order to define a {\it continuous set function} on
a subset of $L^{r(x)}(\Omega)$. We still denote each closed ball centered in the origin of radius $R$ from $L^{r(x)}(\Omega)$
by $\overline B_R(0)$, i.e.
$$\overline B_R(0):=\{u\in L^{r(x)}(\Omega);\;\;|u|_{r(x)}\leq R\}\,.$$
By  Theorem \ref{t3} b) we deduce that $A_\star$ is a continuous function. By the proof of Theorem \ref{t3} a) we have $A_\star(0)>0$
and there exists $R_1>0$ such that $A_\star(R_1)< 0$. Thus, we infer that there exists $R_0>0$ such that $A_\star(R_0)=0$.

We define
$$\Gamma=\{\overline B_R(0)\setminus\overline B_{R_0}(0);\;R\geq R_0\}\subset L^{r(x)}(\Omega)\,$$
and $\mu:\Gamma\rightarrow[0,\infty)$ by
$$\mu(\overline B_R(0)\setminus\overline B_{R_0}(0))=-A_\star(R),\;\;\;\forall\;R\geq R_0\,.$$
By Theorem \ref{t3} a) we find that function $\mu$ has the following properties:
\smallskip

\noindent 1) $\mu(\emptyset)=0$;
\smallskip

\noindent 2) For each $S_1$, $S_2\in \Gamma$ such that $S_1\subset S_2$ we have $\mu(S_1)\leq \mu(S_2)$.

Thus, $\mu$ is a {\it set function} on $\Gamma$. By Theorem \ref{t3} b) and Remark 4 we have that for each $S\subset \Gamma$ it holds true that
$$\mu(S)=\sup_{T\subseteq S}\mu(T)\;\;\;{\rm and}\;\;\;\mu(S)=\inf_{T\supseteq S}\mu(T)\,.$$
We conclude that $\mu$ is a {\it continuous set function} on $\Gamma$.
\medskip

\noindent{\bf Acknowledgments.} The first two authors have been supported by Grant CNCSIS 79/2007 ``Degenerate and Singular Nonlinear Processes".


\begin{thebibliography}{99}
{\footnotesize

%\bibitem{acerbi2} E. Acerbi and G. Mingione, Gradient estimates for the $p(x)$-Laplacean system,
%{\it J.~Reine Angew. Math.} {\bf 584} (2005), 117-148.

\bibitem{AH} M. S. Ashbaugh and E. M. Harrell, Maximal and minimal eigenvalues and their associated nonlinear equations,
{\it J. Math. Phys.} {\bf 28} (1987), 1770-1786.

\bibitem{AHed} D.R. Adams and L.I. Hedberg, {\it
Function Spaces and Potential Theory},
Grundlehren der Mathematischen Wissenschaften [Fundamental Principles of Mathematical Sciences], 314, Springer-Verlag, Berlin, 1996.

\bibitem{A} R. Adams, {\it Sobolev Spaces}, Academic Press, New York,
1975.

\bibitem{bonder} J. F. Bonder and L. M. Del Pezzo, An optimization problem for the first eigenvalue
of the $p$-Laplacian plus a potential, {\it Communication on Pure and Applied Analysis} {\bf 5} (2006),
675-690.

\bibitem{B} H. Brezis, {\it Analyse fonctionnelle: th\'eorie et
applications}, Masson, Paris, 1992.

%\bibitem{CLR} Y. Chen, S. Levine and M. Rao, Variable exponent, linear growth functionals in image
%processing, {\it SIAM J.~Appl. Math.} {\bf 66} (2006), 1383-1406.

\bibitem{Clem1} Ph. Cl\'ement, M. Garc\'{\i}a-Huidobro,
R. Man\'asevich, and K. Schmitt, Mountain pass type solutions for
quasilinear elliptic equations, {\it Calc. Var.} {\bf 11} (2000),
33-62.

\bibitem{Clem2}  Ph. Cl\'ement, B. de Pagter, G. Sweers, and
F. de Th\'elin, Existence of solutions to a semilinear elliptic
system through Orlicz-Sobolev spaces, {\it Mediterr. J. Math.}
{\bf 1} (2004), 241-267.

%\bibitem{D} L. Diening, {\it Theoretical and Numerical Results for
%Electrorheological Fluids}, Ph.D. thesis, University of Frieburg,
%Germany, 2002.

\bibitem{edm} D. E. Edmunds, J. Lang, and A. Nekvinda, On $L^{p(x)}$
norms, {\it Proc. Roy. Soc. London Ser.~A} {\bf 455} (1999), 219-225.

\bibitem{edm2} D. E. Edmunds and J. R\'akosn\'{\i}k, Density of smooth
functions in $W^{k,p(x)}(\Omega)$, {\it Proc. Roy. Soc. London Ser.~A}
{\bf 437} (1992), 229-236.

\bibitem{edm3} D. E. Edmunds and J. R\'akosn\'{\i}k, Sobolev embedding
with variable exponent, {\it Studia Math.} {\bf 143} (2000),
267-293.

\bibitem{egn} H. Egnell, Extremal Properties of the First Eigenvalue of a Class of Elliptic
Eigenvalue Problems, {\it Annali Scuol. Norm. Sup. Pisa} {\bf 14} (1987), 1-48.

\bibitem{F} X. Fan, Remarks on eigenvalue problems involving the $p(x)$-Laplacian,
{\it J. Math. Anal. Appl.} {\bf 352} (2009), 85-98.

\bibitem{FZZ} X. Fan, Q. Zhang and D. Zhao, Eigenvalues of $p(x)$-Laplacian
Dirichlet problem, {\it J. Math. Anal. Appl.} {\bf 302} (2005), 306-317.

\bibitem{Gar} M. Garci\'a-Huidobro, V.K. Le, R. Man\'asevich,
and K. Schmitt, On principal eigenvalues for quasilinear elliptic
differential operators: an Orlicz-Sobolev space setting,
{\it Nonlinear Differential Equations Appl. (NoDEA)} {\bf 6} (1999), 207-225.

\bibitem{G} J.P. Gossez, Nonlinear elliptic boundary value
problems for equations with rapidly (or slowly) increasing
coefficients, {\it Trans. Amer. Math. Soc.} {\bf 190} (1974),
163-205.

%\bibitem{hal} T. C. Halsey, Electrorheological fluids, {\it Science}
%{\bf 258} (1992), 761-766.

%\bibitem{harj} P. Harjulehto, P. H\"ast\"o, M. Koskenoja and S. Varonen,
%Sobolev capacity on the space $W^{1,p(\cdot)}(R^n)$, {\it J. Funct. Spaces Appl.} {\bf 1} (2003), no. 1, 17-33.

\bibitem{KR} O. Kov\'a\v cik and J. R\'akosn\'{\i}k, On spaces
$L^{p(x)}$ and
$W^{1,p(x)}$, {\it Czechoslovak Math. J.} {\bf 41} (1991), 592-618.

\bibitem{RoyalSoc} M. Mih\u ailescu and V. R\u adulescu,
A multiplicity result for a nonlinear degenerate problem arising in
the theory of electrorheological fluids, {\it Proceedings of the Royal Society A: Mathematical, Physical and Engineering Sciences}
{\bf 462} (2006), 2625-2641.

\bibitem{mihradproc} M. Mih\u ailescu and V. R\u adulescu, On a nonhomogeneous quasilinear
eigenvalue problem in Sobolev spaces with variable exponent, {\it
Proc. Amer. Math. Soc.} {\bf 135} (2007), 2929-2937.

\bibitem{jmaa} M. Mih\u ailescu and V. R\u adulescu, Existence and multiplicity of
solutions for quasilinear nonhomogeneous problems: an Orlicz-Sobolev
space setting, {\it Journal of Mathematical Analysis and
Applications} {\bf 330} (2007), Vol. 1, 416-432.

\bibitem{manuscripta} M. Mih\u ailescu and V. R\u adulescu, Continuous spectrum
for a class of nonhomogeneous differential operators, {\it
Manuscripta Mathematica} {\bf 125} (2008), 157-167.

\bibitem{AA} M. Mih\u ailescu and V. R\u adulescu, Eigenvalue problems associated to nonhomogeneous differential
operators in Orlicz-Sobolev spaces, {\it Analysis and Applications} {\bf 6}  (2008), No. 1, 1-16.

\bibitem{fourier} M. Mih\u ailescu and V. R\u adulescu, Neumann
problems associated to nonhomogeneous differential operators in
Orlicz-Sobolev spaces, {\it Ann. Inst. Fourier} {\bf 58} (6)
(2008), 2087-2111.

\bibitem{blms} M. Mih\u ailescu and V. R\u adulescu, Spectrum
in an unbounded interval for a class of nonhomogeneous
differential operators, {\it Bulletin  of the London Mathematical
Society} {\bf 40} (6) (2008), 972-984.

\bibitem{M} J. Musielak, {\it Orlicz Spaces and Modular  Spaces},
Lecture Notes in Mathematics, Vol. 1034, Springer-Verlag, Berlin,
1983.

\bibitem{rao} M.M. Rao and Z.D. Ren, {\it Theory of Orlicz Spaces}, Marcel Dekker, Inc.,
New York, 1991.

%\bibitem{R} M. Ruzicka, {\it Electrorheological Fluids: Modeling
%and Mathematical Theory}, Springer-Verlag, Berlin, 2002.

\bibitem{samko} S. Samko and B. Vakulov, Weighted Sobolev theorem with variable exponent
for spatial and spherical potential operators, {\it J. Math. Anal. Appl.} {\bf 310} (2005), 229-246.

\bibitem{S} M. Struwe, {\it Variational Methods: Applications to
Nonlinear Partial Differential Equations and Hamiltonian Systems},
Springer, Heidelberg, 1996.

%\bibitem{Z1} V. Zhikov, Averaging of functionals in the calculus of
%variations and elasticity, {\it Math. USSR Izv.} {\bf 29} (1987),
%33-66.

}
\end{thebibliography}
\end{document}